\documentclass{article}
\pdfoutput=1
\usepackage{amssymb, latexsym, amsmath, verbatim, amsthm, amscd}
\usepackage{pinlabel}
\usepackage{graphicx}
\usepackage[pdftex]{hyperref}

\def\mytitle{Computing the Maximum Slope Invariant in Tubular Groups}

\title{\mytitle}
\author{Christopher H. Cashen\\
Department of Mathematics\\
University of Utah\\
Salt Lake City, UT 84112\\
\href{mailto:cashen@math.utah.edu}{\nolinkurl{cashen@math.utah.edu}}\\
\href{http://www.math.utah.edu/~cashen}{\nolinkurl{http://www.math.utah.edu/~cashen}}}

\date{\today}

\hypersetup{
    pdftitle={Computing the Maximum Slope Invariant in Tubular Groups},    
    pdfauthor={Christopher H. Cashen},     
    pdfkeywords={Tubular Groups}, 
    colorlinks=true,       
    linkcolor=black,          
    citecolor=black,        
    filecolor=black,      
    urlcolor=black           
}

\theoremstyle{plain}

\theoremstyle{remark}

\theoremstyle{definition}

\def\makeautorefname#1#2{\expandafter\def\csname#1autorefname\endcsname{#2}}
\let\fullref\autoref

\makeautorefname{theorem}{Theorem} 
\makeautorefname{lemma}{Lemma} 
\makeautorefname{proposition}{Proposition} 
\makeautorefname{corollary}{Corollary} 
\makeautorefname{definition}{Definition}
\makeautorefname{example}{Example}
\makeautorefname{subsection}{Section}
\makeautorefname{subsubsection}{Section}
\makeautorefname{section}{Section}

\makeatletter 
\let\c@lemma=\c@theorem 
\makeatother
\makeatletter 
\let\c@proposition=\c@theorem 
\makeatother
\makeatletter 
\let\c@corollary=\c@theorem 
\makeatother
\makeatletter 
\let\c@definition=\c@theorem 
\makeatother
\makeatletter 
\let\c@example=\c@theorem 
\makeatother

\DeclareMathOperator{\slope}{slope}

\DeclareMathOperator{\Z}{\mathbb{Z}}

\def\from{\colon\thinspace}

\newcommand{\case}[2]{\smallskip \noindent \emph{#1}:\newline #2 \medskip}

\begin{document}

\maketitle

\begin{abstract}
  We show that the maximum slope invariant for tubular groups is easy to
  calculate, and give an example of two tubular groups that are
  distinguishable by their maximum slopes but not by edge pattern
  considerations or isoperimetric function.
\end{abstract}



\section{Introduction}
The main examples of tubular groups are constructed from 
$\mathbb{Z}^2$ by amalgamating along cyclic subgroups.

Tubular groups have been used as examples of various phenomena. Brady and
Bridson computed isoperimetric functions for the groups
\[BB(p,r)=\left<a,b,x,y\mid [a,b]=1,\, x^{-1}a^px=a^rb,\,y^{-1}a^py=a^rb^{-1}\right>\]
with $0<p<r$ and found that the degrees of the isoperimetric functions form a dense set in
$[2,\infty)$. 
The Dehn function of $BB(p,r)$ is $x^n$ with $n=2\log_2\frac{2r}{p}$ \cite{BraBri00}.

Later Brady, Bridson, Forester and Shankar, \cite{BraBriFor09}, used
more complicated examples of tubular groups, called ``snowflake groups'', to show that  all the rationals in that
interval appear. 

Other examples include the examples of Croke and Kleiner \cite{CroKle00}
and Wise's group \cite{Wis96}. 

There are two quasi-isometry invariants readily available among these groups. For the Brady-Bridson and snowflake groups we have the isoperimetric exponent. Mosher, Sageev and Whyte give another: the collection of affine equivalence classes of edge patterns \cite{MosSagWhy04}.

In a previous paper \cite{Cas07b} we gave an algorithm to decide
whether or not two tubular groups are quasi-isometric. In the
course of the proof we constructed a tree, the ``tree of P--sets'',
related to the Bass-Serre tree of a particular graph of groups
decomposition of the group. 
The (directed) edges of this tree come with a parameter, the
``height change across the edge''. The maximum slope invariant for the tubular
group is then the maximum coarse slope of any ray in the tree, where
coarse slope is average height change per unit length, in an
appropriate sense. It follows from the quasi-isometry algorithm that
the maximum coarse slope is a quasi-isometry invariant of the group
(at least in certain cases).

Among the groups $BB(p,r)$ the edge patterns are all equivalent and it turns out that the isoperimetric exponent is a complete quasi-isometry invariant. 
One might wonder whether the equivalence classes of edge patterns and the isoperimetric exponent always determine the quasi-isometry class of a tubular group. 

The answer is ``no''. In this note we give an example of two tubular groups with the same equivalence classes of edge patterns and the same isoperimetric exponent but different maximum slopes.

The method of computing the maximum slopes is easy and generalizes to many classes of tubular groups.
Dehn functions for tubular groups are known only for special cases,
and  running the quasi-isometry algorithm of \cite{Cas07b} is
laborious, so, in practice, checking that maximum slopes are equal is a
good test to perform before running the quasi-isometry algorithm. 

\section{Preliminaries}
In this section we recall the necessary terminology and facts about
tubular groups and quasi-isometries of tubular groups.

We simplify the exposition considerably by considering
only tubular groups that are graphs of groups with
vertex groups $\mathbb{Z}^2$ and edge groups
$\mathbb{Z}$ such that the edge groups incident to a vertex inject into exactly three distinct maximal cyclic subgroups.

The reader is referred to \cite{Cas07b} for more background in tubular
groups and graphs of groups, and for proofs of claims in this section.

\subsection{The Geometric Model}
Consider a graph of groups with
vertex groups $\mathbb{Z}^2$ and edge groups
$\mathbb{Z}$  such that the edge groups incident to a vertex inject into exactly three distinct maximal cyclic subgroups. 
A model for such a group can be built by taking a torus
for each vertex, an annulus for each edge, and gluing a boundary of an
annulus to a torus according to the corresponding edge injection.

If we then choose a metric on each torus and annulus and lift to the
universal cover we get a geodesic metric space quasi-isometric to the
group, which we call the \emph{geometric model}.

Consider one of the vertex $\Z^2$ groups. 
Suppose $\Z^2=\left<a,b\mid[a,b]=1\right>$.
The ``usual'' metric is the one for which the word $a^xb^y$ in
$\mathbb{Z}^2$ corresponds to the point with Cartesian coordinates
$(x,y)$ in the standard Euclidean plane. 
Suppose the incident edge groups inject into the maximal cyclic subgroups containing $a$, $a^rb^s$, and $a^tb^u$, where $0$, $\frac{s}{r}$ and $\frac{u}{t}$ are distinct. 
This can be arranged by picking a new generating set, if necessary, since we have assumed that there are three distinct maximal cyclic subgroups.

This means that with this choice of coordinates, in the universal cover the
edge spaces incident to this vertex space attach along lines of slope
$0$, $\frac{s}{r}$ and $\frac{u}{t}$.
For any of these slopes, since the edge group $\mathbb{Z}$ is
infinite index in the vertex group $\mathbb{Z}^2$, in the universal
cover we have infinitely many edge spaces gluing on to the vertex
space along lines of the same slope.

We will have three different infinite families of parallel lines in
the plane. This collection of families of infinite lines is called the
\emph{edge pattern} in the vertex space.

From work of Mosher, Sageev and Whyte \cite{MosSagWhy04} we know that
a quasi-isometry of tubular groups takes vertex groups to within
bounded distance of vertex groups. We also know that these edge
patterns are the coarse intersection patterns of the various vertex
spaces, and must be preserved up to bounded error by
quasi-isometries.
Furthermore, when there are infinite families of at least three
different slopes in a plane, then only quasi-isometries of the plane
that preserve the three families are bounded distance from an affine
map.
Therefore, the affine equivalence class of the set of slopes in a
vertex space is a quasi-isometry invariant of the group. 

For our examples we may assume that any quasi-isometry
restricts on each vertex space to an affine map that takes three
specified slopes to three specified slopes. 
Projectively there is a unique map that does this, so the only freedom
in the map is translation and rescaling the entire plane by a constant.
Moreover, there is a convenient choice of metric that will make
quasi-isometries on the vertex spaces particularly nice.

Choose the metric on the vertex group so that the word $a^xb^y$ corresponds to the vector $A(_y^x)$ in the plane, where $A$ is the matrix:
\[A=\left(\begin{matrix}
1 & -\frac{1}{2}\frac{ru+st}{su}\\
0 & \frac{\sqrt{3}}{2}\frac{ru-st}{su}
\end{matrix}\right)\]

This is a convenient choice because it makes the line pattern
symmetric, the three families of parallel lines differ from one
another by angle $\frac{\pi}{3}$.
Any permutation of the slopes can be achieved by an
isometry of the plane.

Once the metrics have been chosen on the vertex groups we can define
height change across an edge. 
Each edge group is $\mathbb{Z}$. 
Define the stretch factor across the (directed) edge to be the ratio
of the lengths of the image of the generator of the edge group in the
two adjacent vertex groups.
The \emph{height change} across the edge is $-\log_2(\text{stretch factor})$.

\subsection{Quasi-isometries preserve height change}
For $i=1,2$, let $G_i$ be a tubular group, $D_i$ its Bass-Serre tree, and $X_i$ its
geometric model. Let $q_i \from X_i \to D_i$ be the usual quotient map.
A quasi-isometry $\phi\from X_1\to X_2$ induces a bijection $\phi_*$
from vertices of $D_1$ to vertices of $D_2$.
This bijection of vertices can be extended to a continuous map $D_1\to
D_2$ by connecting the dots.

The height change between two vertices of the Bass-Serre tree of a
tubular group is the sum of the height changes across the edges of the
geodesic segment joining the vertices.

The quasi-isometry $\phi$ is coarsely height preserving in the sense
that there is some constant $C$ so that for any two vertices $v, v'\in
D_1$, the height change between $\phi_*(v)$ and $\phi_*(v')$ is within
$C$ of the height change between $v$ and $v'$.

\subsection{Coarse Slope}
Let $c\from [0,\infty)\to D_1$ be a geodesic ray.
A quasi-isometry $\phi\from X_1\to X_2$ induces a bijection $\phi_*$
from vertices of $D_1$ to vertices of $D_2$, but does not necessarily
preserve adjacency.
Thus, the ray $\phi_*\circ c\from [0,\infty)\to D_2$ may not be a geodesic ray.

Let $e$ and $e'$ be edges incident to a vertex $v$ in $D_i$. 
We say that $e$ and $e'$ are \emph{parallel at $v$} if
$q_i^{-1}(e)\cap q_i^{-1}(v)$ and $q_i^{-1}(e')\cap q_i^{-1}(v)$ are
parallel lines in $q_i^{-1}(v)$. 
Equivalently, $e$ and $e'$ are parallel at $v$ if their stabilizer
subgroups in $G_i$ are contained in a common maximal cyclic subgroup
of the stabilizer of $v$.

We say that $c$ \emph{has a twist at $t$} (or at $v=c(t)$) if the incoming
and outgoing edges at $c(t)$ are not parallel at $c(t)$.

Quasi-isometries preserve twists: $c$ has a twist at $t$ if and only
if $\phi_*\circ c$ has a twist at $t$. 

So, while the length of a segment in $D_1$ is not preserved by
$\phi_*$, the number of twists along that segment is, and we define
the coarse slope of a ray $c$ in $D_1$ to be the ratio of height
change to number of twists.

Let $\mathrm{Twist}_c(t)$ denote the number of twists of $c(0,t)$.
A ray $c$ in $D_i$ \emph{has coarse slope} $m$, $\slope (c)=m$, if
there exists a $C>0$ such that for all $t$, \[m\cdot \mathrm{Twist}_c(t)-C\leq
h(c(0),c(t)) \leq m\cdot \mathrm{Twist}_c(t)+C\]

Not every ray has a coarse slope, but $c$ has coarse slope $m$ if and
only if $\phi_*\circ c$ does.
The \emph{maximum slope} invariant of a tubular group is the
quasi-isometry invariant given by the largest number
$m$ (possibly $\infty$) that occurs as the coarse slope of a ray in the Bass-Serre tree.
(The supremum of coarse slopes is always achieved because the group
acts cocompactly on the Bass-Serre tree.)

\section{The Procedure}\label{sec:procedure}
In this section we give the procedure for computing the maximum slope, using the group $BB(p,r)$ to illustrate.

\case{Step 1: Compute the Height Changes}{
Start with a graph of groups decomposition. Call the vertices $v_1,\dots v_j$. 
Choose the metric that makes the edge patterns symmetric and compute height changes across each edge.

\begin{figure}[h]
\labellist
\large
\pinlabel{$\Z^2$} at 61 780
\pinlabel{$a^p$} at 53 800
\pinlabel{$a^p$}   at 69 800
\pinlabel{$a^rb$}   at 53 757.8
\pinlabel{$a^rb^{-1}$}   at 69 758
\endlabellist
 \centering
 \includegraphics[width=.5\textwidth]{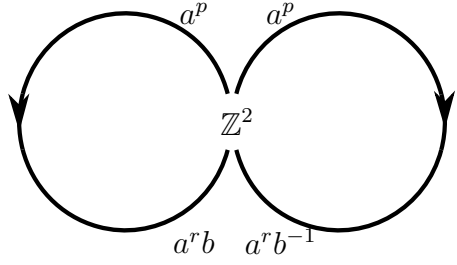}
 \caption{graph of groups for $BB(p,r)$}
 \label{fig:bb}
\end{figure}

The graph of groups for $BB(p,r)$ is depicted in
\fullref{fig:bb}. Let $x$ and $y$ be stable letters corresponding to the loop on the left and right, respectively.
The labels at the ends of the edges indicate that these elements of
the vertex group are conjugate by the stable letter associated to the
edge, thus $x$ conjugates $a^p$ to $a^rb$.
The fundamental group of the graph of groups is:
\[BB(p,r)=\left<a,b,x,y\mid [a,b]=1,\, x^{-1}a^px=a^rb,\,y^{-1}a^py=a^rb^{-1}\right>\]

The matrix $A$ that determines the symmetric metric on the vertex group is:
\[A=\left(\begin{matrix}
1 &0\\
0 & r\sqrt{3}
\end{matrix}\right)\]

The height change across $x$ is
\[H=-\log_2\frac{|A(^r_1)|}{|A(^p_0)|}=-\log_2\frac{2r}{p}\]
The height change across $y$ is also $H$. Note that $H$ is negative since $r>p$.
}

\case{Step 2: Identify Parallel Edges}{
At each vertex, join edges with an arc if their groups inject into a
common maximal cyclic subgroup.
This happens exactly when the edges are parallel at the vertex.

If there is a loop in the graph such that for every edge-vertex-edge sequence in the loop the two edge groups are parallel and such that the net height change around the loop is not zero, stop. The maximum slope is infinite.

\begin{figure}[h]
\labellist
\large
\pinlabel{$H$}   at 15 781
\pinlabel{$H$}   at 105 781
\endlabellist
\centering
\includegraphics[width=.5\textwidth]{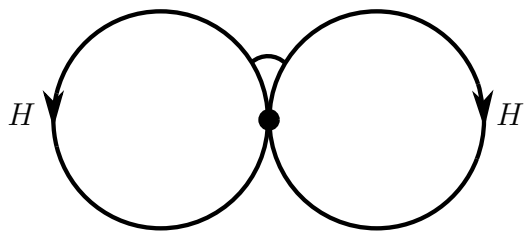}
\caption{height change and parallel edges for $BB(p,r)$}
\label{fig:bb2}
\end{figure}

In $BB(p,r)$ the two edge groups both map into the maximal cyclic
subgroup $\left<a\right>$ at the top of \fullref{fig:bb2}. 
There are no non-trivial loops for which all edge-vertex-edge transitions involve parallel edges.
}

\case{Step 3: Fold Parallel Edges}{
Subdivide each edge by adding a vertex at the midpoint. Call these new vertices $m_1,\dots m_k$.

\begin{figure}[h]
\centering
\includegraphics[width=.5\textwidth]{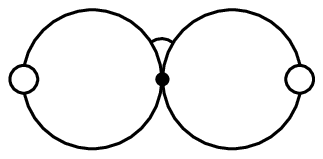}
\caption{subdivided edges for $BB(p,r)$}
\label{fig:bb3}
\end{figure}

At each $v_i$, for each collection of edges incident to the vertex and joined by arcs, fold them together by identifying them up to their midpoints. 

In the resulting \emph{graph of P--sets}, label each edge with a height change in such a way that the height changes between the original vertices is preserved. 
This is always possible.
One way to accomplish this is to look at each remaining $m_i$. 
It is adjacent to some of the original vertices $v_{i_1},\dots, v_{i_l}$. 
Among these there is some $\alpha$ so that for all $\beta$, the height change from $v_{i_\alpha}$ to $v_{i_\beta}$ is non-negative. 
Give the edge from $v_{i_\alpha}$ to $m_i$ height change 0, and give the edge from $m_i$ to $v_{i_\beta}$ the height change equal to the height change from $v_{i_\alpha}$ to $v_{i_\beta}$. 

One may check that this graph of P--sets is the quotient by the group
action on the tree of P--sets constructed in \cite{Cas07b}.

We will leave edges with height change 0 unlabeled and undirected.

\begin{figure}[h]
\labellist
\large
\pinlabel{$|H|$}   at 68 800
\endlabellist
\centering
\includegraphics[width=.5\textwidth]{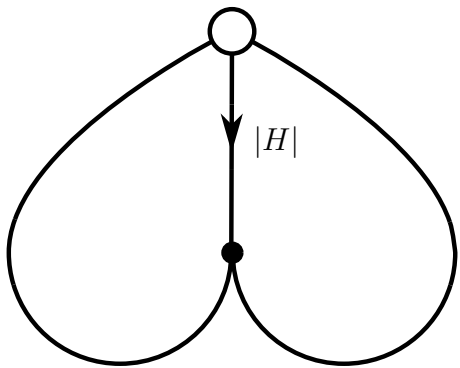}
\caption{graph of P--sets for $BB(p,r)$}
\label{fig:bb4}
\end{figure}

}

\case{Step 4: Find the Embedded Loop of Maximum Slope}{
The effect of this procedure is to collapse parallel edges at a vertex
to a single edge. 
Since we assumed that for each vertex there are exactly three maximal
cyclic subgroups containing the edge injections, the original (small
black) vertices all have valence 3. 
The new (large white) vertices represent collections of parallel
edges. 

Let $c$ be a path in the graph of P--sets and let $v$ be an original vertex in
the interior of the path.
For any lift $\tilde{c}$ of $c$ to the Bass-Serre tree, $\tilde{c}$
has a twist at the vertex $\tilde{v}$ corresponding to $v$.
This means that length of a path in the graph of P--sets corresponds to number of
twists in any lift of the path, so we can compute slopes by taking the
ratio of height change to length of paths in the graph.
(Here we consider edges in the graph to have length one half.)

If the graph of P--sets is a tree then the maximum slope is 0. Otherwise, the maximum slope can be realized by an embedded loop. 
To see this, consider a loop that is not embedded and has non-zero slope. 
It contains some embedded sub-loops, and the slope of the entire loop is less than or equal to the maximum of the slopes of the embedded sub-loops. 

For $BB(p,r)$ the maximum slope is $|H|=\log_2\frac{2r}{p}$.
Recall that the isoperimetric exponent for $BB(p,r)$ is
$2\log_2\frac{2r}{p}$, which is exactly twice the maximum slope.
}

\section{An Example}
Consider the snowflake group, $G=G_{r,P}$ of \cite{BraBriFor09}, 
with $P=(4)$, $p=4$, $q=1$. Let $r=2^p=16$.
This snowflake group has Dehn function $x^4$.

$G$ is the fundamental group of a graph of groups with three
$\mathbb{Z}^2$ vertex groups
$\mathbb{Z}^2=\left<a_i,b_i\mid [a_i,b_i]=1\right>$ for $i=1,2,3$.

Each vertex has edges injecting into the cyclic subgroups generated by $a_i$, $b_i$, and $a_ib_i$.

\begin{figure}[h]
\labellist
\large
\pinlabel{$\Z^2$}   at 66 731
\pinlabel{$\Z^2$}   at 106 731
\pinlabel{$\Z^2$}   at 146 731
\pinlabel{$c$}   at 55 725
\pinlabel{$b_1$}   at 76 725
\pinlabel{$a_2b_2$}   at 92 725
\pinlabel{$b_2$}   at 115 725
\pinlabel{$a_3b_3$}   at 131 725
\pinlabel{$a_4^{16}$}   at 157 725
\pinlabel{$a_1^{16}$} [l] at 66 740
\pinlabel{$a_2^{16}$} [l] at 106 740
\pinlabel{$a_3^{16}$} [l] at 146 740

\endlabellist
  \centering
  \includegraphics[width=.8\textwidth]{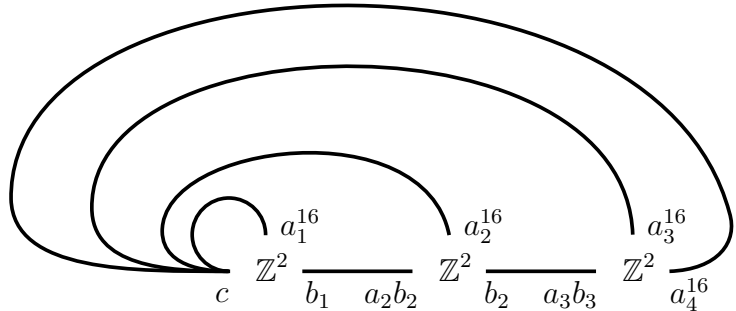}
  \caption{A graph of groups diagram for the snowflake group $G_{16, (4)}$}
  \label{Fig:snowflake}
\end{figure}

\fullref{Fig:snowflake} gives a graph of groups diagram for this
group.
In the diagram, $a_4=b_3$ and $c=a_1b_1=a_1a_2a_3a_4$.

Conveniently, symmetrization in this case does not stretch the lines,
so stretch factors are just ratios of indices. 
The edge that goes from $a_i^{16}$ to $c$ therefore has stretch factor $\frac{1}{16}$, and height change $-\log_2\frac{1}{16}=4$.
The other edges have 0 height change.

\begin{figure}[h]
\labellist
\large
\pinlabel{4} [l] at 63 746
\pinlabel{4} [b] at 75 757
\pinlabel{4} [b] at 75 776
\pinlabel{4} [b] at 75 790

\endlabellist
  \centering
  \includegraphics[width=.8\textwidth]{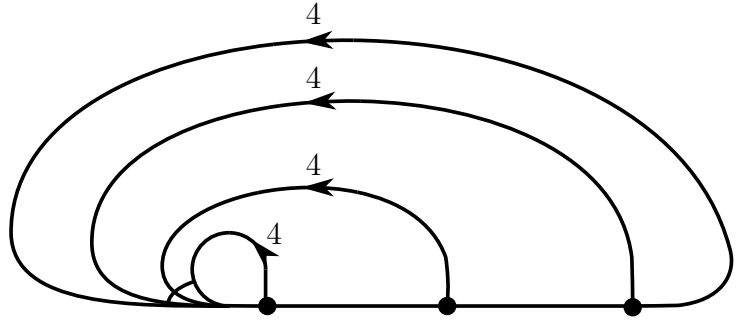}
  \caption{Height changes and parallel edges for the snowflake group $G_{16, (4)}$}
  \label{Fig:snowflakewheights}
\end{figure}

In \fullref{Fig:snowflakewheights} we have height changes and parallel
edges for $G_{16, (4)}$.

\begin{figure}[h]
\labellist
\large
\pinlabel{$4$}   at 52 726
\endlabellist
  \centering
  \includegraphics[width=.8\textwidth]{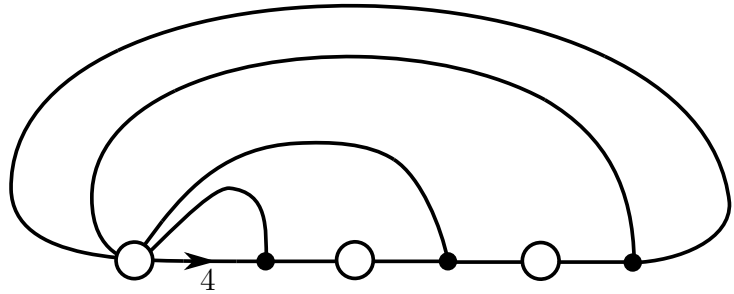}
  \caption{The P--set graph for $G_{16,(4)}$}
  \label{fig:snowflakePgraph}
\end{figure}

Now if we fold this graph down to a graph of P--sets, \fullref{fig:snowflakePgraph}, we quickly see that the maximum slope in the tree of P--sets for this group is 4 (the smallest loop in the figure realizes the maximum slope).
The isoperimetric exponent was also 4.

However, recall that $BB(p,r)$ has isoperimetric exponent equal to twice the maximum slope. $BB(1,2)$ has Dehn function $x^4$ but a maximum slope of 2. It is not quasi-isometric to the snowflake group.

\section{Generaliztion}
The procedure described in \fullref{sec:procedure} generalizes to
arbitrary tubular groups, subject to constraints described in
\cite{Cas07b}. 
First, coarse slope is only well defined when the edge patterns in
every vertex space have at least three distinct families of parallel lines. 
Second, the slopes depend on the choice of symmetric metric on the
vertex spaces. 
For vertex spaces whose edge patterns consist of three
or four slopes there is a unique choice, up to isometry and rescaling. 
For an edge pattern with five or more families of parallel lines this
is not always true, there may or may not be a
unique choice of symmetric metric.
The maximum coarse slope is still useful in distinguishing groups when
there are edge patterns without a canonical choice of symmetric
metric, but one must take care to identify each pattern in an
equivalence class with the same symmetric pattern.



\end{document}